%% file: root.tex
\documentclass[letterpaper, 10 pt, conference]{ieeeconf}  

\input{extras/common_commands.tex}

\input{extras/common_preamble.tex}
\input{extras/notation.tex}

\IEEEoverridecommandlockouts                              

\usepackage{graphics} 

\title{\LARGE \bf
  Robust Model Predictive Control for Linear Systems with State
  and Input Dependent Uncertainties
}

\author{Danylo Malyuta$^{1}$, \Behcet{} \Acikmese{}$^{2}$ and Martin Cacan$^{3}$%
  \thanks{$^{1}$Ph.D. student, W.E. Boeing Department of Aeronautics \&
    Astronautics, University of Washington, Seattle, WA 98195, USA
    \texttt{danylo@uw.edu}}%
  \thanks{$^{2}$Professor, W.E. Boeing Department of
    Aeronautics \& Astronautics, University of Washington, Seattle, WA 98195,
    USA \texttt{behcet@uw.edu}}%
  \thanks{$^{3}$Guidance and Control Analyst, Jet Propulsion Laboratory, California Institute of
    Technology, Pasadena, CA 91109, USA \texttt{martin.cacan@jpl.nasa.gov}}%
}

\begin{document}

\maketitle
\thispagestyle{empty}
\pagestyle{empty}

\begin{abstract}
  This paper presents a computationally efficient robust model predictive
  control law for discrete linear time invariant systems subject to additive
  disturbances that may depend on the state and/or input norms. Despite the
  dependency being non-convex, we are able to handle it as a second-order cone
  program. Both open-loop and semi-feedback planning strategies are
  presented. The formulation has linear complexity in the planning horizon
  length. The approach is thus amenable to efficient real-time implementation
  with a guarantee on recursive feasibility and global optimality. Robust
  position control of a satellite is considered as an illustrative example.
\end{abstract}

\input{introduction}
\input{problem_formulation}
\input{control_law}
\input{main_results}
\input{extensions}
\input{example}
\input{conclusion}

\section{Acknowledgment}

This research has been partly supported by the National Science Foundation,
grant number CMMI-1613235, and was partially carried out at the Jet Propulsion
Laboratory, California Institute of Technology, under a contract with the
National Aeronautics and Space Administration. Government sponsorship
acknowledged. The authors would like to extend special gratitude to Daniel
P. Scharf, David S. Bayard, Jack Aldrich and Carl Seubert for their helpful
insight and discussions.

\bibliographystyle{ieeetr}
\bibliography{bibliography.bib}

\end{document}

%% file: extras/common_commands.tex
\usepackage{xcolor}
\definecolor{commentgreen}{RGB}{0, 153, 0}



%% file: extras/common_preamble.tex
\usepackage{amsthm}
\usepackage{amsmath}               
{
  \theoremstyle{plain}
  
  \newtheorem{definition}{Definition}
  \newtheorem{theorem}{Theorem}
  \newtheorem{corollary}{Corollary}
}
\usepackage{amssymb}
\usepackage{bm}
\usepackage{mathtools}
\usepackage{graphicx}
\usepackage{tabularx}
\usepackage[]{hyperref}
\usepackage{xcolor}
\usepackage{algorithm, algorithmicx}
\usepackage[noend]{algpseudocode}

\newlength{\captionskip}


%% file: extras/notation.tex
\newcommand{\BEAS}{\begin{eqnarray*}}
\newcommand{\EEAS}{\end{eqnarray*}}
\newcommand{\BEA}{\begin{eqnarray}}
\newcommand{\EEA}{\end{eqnarray}}
\newcommand{\BEQ}{\begin{equation}}
\newcommand{\EEQ}{\end{equation}}
\newcommand{\BIT}{\begin{itemize}}
\newcommand{\EIT}{\end{itemize}}
\newcommand{\BNUM}{\begin{enumerate}}
\newcommand{\ENUM}{\end{enumerate}}

\newcommand{\BA}{\begin{array}}
\newcommand{\EA}{\end{array}}

\newcommand{\reals}{\mathbb R}

\newcommand{\integers}{\mathbb Z}

\newcommand{\Co}{{\mathop \mathrm{co}}}

{\begin{quote}}{\end{quote}}

\newcounter{oursection}

\newcommand{\minoptconstrained}[3]{
  \begin{aligned}
    & \underset{#1}{\text{minimize}}
    && #2 \\
    & \text{subject to}
    && #3
  \end{aligned}%
}

\definecolor{darkolivegreen}{rgb}{0.33, 0.42, 0.18}

\newcommand{\Behcet}{Beh\c{c}et}
\newcommand{\Acikmese}{A\c{c}{\i}kme\c{s}e}

\newcommand{\transp}{{\scriptscriptstyle\mathsf{T}}}


%% file: introduction.tex
\section{Introduction}
\label{}

In this paper we develop a novel convex formulation for robust model predictive
control (RMPC) of discrete linear time invariant systems with additive state
and/or input dependent uncertainty. 
%
A major advantage of RMPC is its ability to guarantee by design that input and
state constraints are satisfied for all uncertainty realizations. Several
extensive survey papers cover available modeling assumptions and solution
methods \cite{Bemporad2007,Mayne2014}. We restrict our attention to discrete
linear time invariant systems and focus on developing a real-time implementable
algorithm on computationally constrained hardware. We assume a perfect dynamics
model affected by additive bounded uncertainty and present open-loop and
semi-feedback planning strategies \cite[Section~8]{Bemporad2007}. Note that
these strategies merely refer to the method of handling uncertainty within the
RMPC optimization problem. In each case, the on-line implementation is done in
the traditional feedback sense of re-solving and applying the first of the
optimal control inputs at every time step.
Several authors have considered closed-loop formulations
\cite{Lee1997,Bemporad1998,Scokaert1998,Kerrigan2004,Lofberg2003}. However,
these suffer from combinatorial complexity in the problem dimension unless
certain restrictive assumptions are made in terms of cost norm or feedback type.

Closest to our work \cite{Blanchini1990b,Blanchini1990a,Blanchini1990c} and
\cite{Blanchini1994} use pre-computed constraint tightening factors to guarantee
robustness to worst-case uncertainty through a set of linear constraints. The
computational cost is marginally higher than that of nominal MPC and the online
problem is at most a quadratic program (QP). A similar idea is exploited for
nonlinear systems in \cite{Marruedo2002}. Dependent uncertainty has received
some attention in nonlinear MPC \cite{Pin2009,Lee2002,Limon2004}. These
formulations, however, are conservative due to their use of a Lipschitz constant
for constraint tightening, which considers uncertainty magnitude but not
direction The incrementally conic uncertainty model in \cite{Acikmese2010} is
used to assure robustness via a static feedback component obtained by an
off-line linear matrix inequality procedure and an on-line nominal MPC law, as
in tube MPC \cite{Rakovic2009}.

Our main contribution is to present how a state and input dependent uncertainty
model can be included in RMPC while retaining low computational complexity. We
present an open-loop formulation first and, because it can be overly
conservative, extend it to a semi-feedback formulation where a static feedback
gain is embedded into the planning task \cite{Bemporad1998,Cannon2005}. Our
model is a subset of \cite{Acikmese2010} but has the advantage of using the more
computationally efficient second-order cone programming (SOCP) for a robust
solution. Furthermore, unlike \cite{Pin2009,Lee2002,Limon2004}, our method
captures disturbance directionality effects via H\"older's inequality and is
therefore less conservative. Finally, we illustrate how the uncertainty model
can describe the very popular Gates thruster execution-error model for
satellites \cite{Gates1963,Wagner2014,Goodson2013}. Because our approach has
linear complexity in the horizon length and is at most an SOCP problem, we
expect the algorithm to be amenable to real-time implementation
\cite{Dueri2014,Dueri2017,Zeilinger2014,Domahidi2013}.


The paper is organized as follows. In Section~\ref{sec:problem_formulation} the
uncertainty model and the resulting robust optimal control problem are
introduced. In Section~\ref{sec:control_law} a solution is presented as an
open-loop RMPC law that is at most an SOCP
problem. Section~\ref{sec:main_results} proves recursive feasibility and
suggests a computationally efficient check of robust controlled
invariance. Section~\ref{sec:extensions} converts the formulation to a
semi-feedback implementation. An example is presented in
Section~\ref{sec:example} which illustrates the method's
effectiveness. Section~\ref{sec:conclusion} discusses possible extensions and
offers concluding remarks.

\textit{Notation}: $\reals$, $\reals_+$ and $\reals_{++}$ are the reals,
non-negative reals and positive reals. Unless otherwise specified, matrices are
uppercase (e.g. $A$), scalars and vectors are lowercase (e.g. $x$) and sets are
calligraphic uppercase (e.g. $\mathcal S$).
$0_{n\times m}\in\reals^{n\times m}$ and $I_n\in\reals^{n\times n}$ are the zero
and identity matrices respectively, where $n=m=3$ when the subscripts are
omitted. Parentheses denote vertical stacking,
e.g. $(1,2,3)\in\reals^3$. $M_i^\transp$ denotes the $i$-th row of matrix
$M$. $\|\cdot\|_p$ denotes the $p$-norm, e.g. $p=1,2,\infty$. $\Co\mathcal S$ is
the convex hull of $\mathcal S$. Direct set mapping is written as
$M\mathcal S=\{Mx : x\in\mathcal S\}$. The support function of
$\mathcal S\subseteq\reals^n$ along a direction $v\in\reals^n$ is denoted
$\sigma_{\mathcal S}(v)\triangleq\max_{z\in\mathcal S}v^\transp z$. The
shorthand $a:b$ represents the integer sequence $a,\dots,b$. $\mathcal S^a$ is
the Cartesian product performed $a$ times,
e.g. $\mathcal S^3 = \mathcal S\times\mathcal S\times\mathcal S$.


%% file: problem_formulation.tex
\section{Problem Formulation}
\label{sec:problem_formulation}

\begin{figure}
  \centering
  \includegraphics[width=0.7\columnwidth]{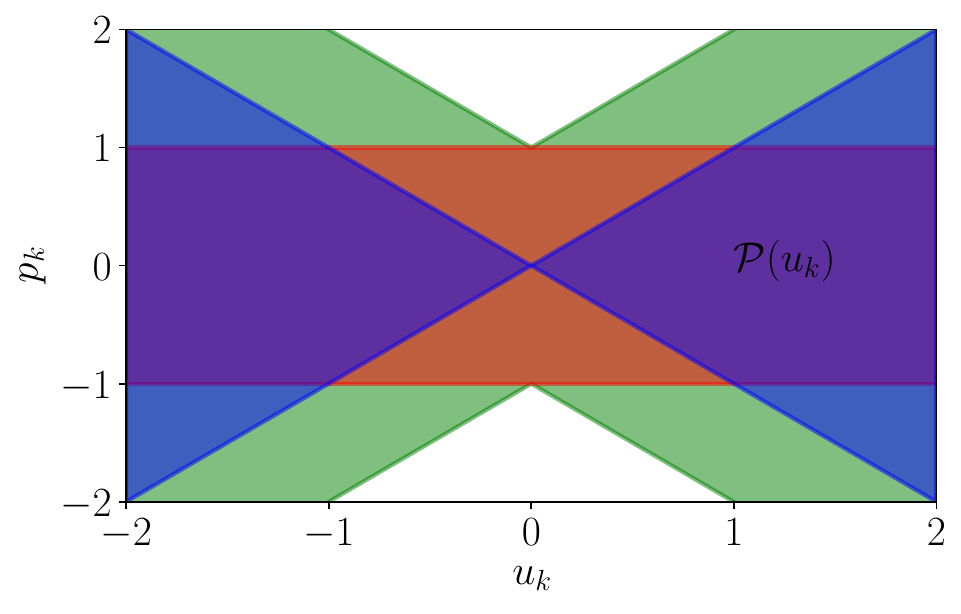}
  \vspace{\captionskip}
  \caption{Non-convex input dependent uncertainty set $\mathcal P(u_k)$ (green),
    expressed as the Minkowski sum of a polytopic independent component (red)
    and a conic dependent component (blue).}
  \label{fig:uncertainty_set_example}
\end{figure}

This section describes the control problem and, in particular, defines the state
and input dependent uncertainty model.  Consider a discrete linear time
invariant system with additive uncertainty:
\begin{equation}
  \label{eq:DLTV}
  x_{k+1} = Ax_k+Bu_k+Dp_k,
\end{equation}
where $k\in\integers_+$, $x_k\in\reals^n$, $u_k\in\reals^m$ and $p_k\in\reals^d$
are, respectively, the time step, state, input and uncertainty, while $A$, $B$
and $D$ are constant matrices of commensurate dimension. The following state and
input constraints are to be respected for all time:
\begin{subequations}
  \begin{alignat}{3}
    x_k &\in \mathcal X \subset \reals^n\quad &&\forall k\in\integers_+,
    \label{eq:state_constraint} \\
    u_k &\in \mathcal U \subset \reals^m\quad &&\forall k\in\integers_+,
    \label{eq:input_constraint} 
  \end{alignat}
\end{subequations}
where $\mathcal X$ and $\mathcal U$ are compact convex polytopes that we express
as follows:
\begin{subequations}
  \begin{align}
    \mathcal X &\triangleq \{x\in\reals^n : Fx\le f\}, \\
    \mathcal U &\triangleq \{u\in\reals^m : Hu\le h\},
  \end{align}
\end{subequations}
where the rows of $F$, $H$ and the elements of $f$, $h$ define respectively the
polytope facet normals and distances. Consider now the following uncertainty
model:
\begin{subequations}
  \begin{alignat}{3}
    p_k
    &\in \mathcal P&&(x_k,u_k)\subset\reals^d, \label{eq:uncertainty_constraint} \\
    \mathcal P(x_k,u_k)
    &\triangleq \{&&Ww+\sum_{i=1}^{n_q}L_iq_i\in\reals^d :
    Rw\le r, \label{eq:uncertainty_set} \\
    & &&\|q_i\|_{p_{q,i}}\le \phi_i(\|F_{x,i}x_k\|_{p_{x,i}},\|F_{u,i}u_k\|_{p_{u,i}}), \nonumber \\
    & &&i=1,\dots,n_q\}, \nonumber
  \end{alignat}
\end{subequations}
where each function $\phi_i:\reals_+\times\reals_+\to\reals_+$ is convex and
non-decreasing such that $\phi_i(\|\cdot\|_{p_{x,i}},\|\cdot\|_{p_{u,i}})$ is
convex \cite{Boyd2004}. The uncertainty set is effectively the sum of an
independent component originating from a polytope and a state and input
dependent component that is bounded by a non-convex inequality. Because
unbounded disturbances are not practical, we assume that $\mathcal P(x_k,u_k)$
is compact $\forall x_k\in\mathcal X,~u_k\in\mathcal U$. Note that
$\mathcal P(x_k,u_k)$ is a non-convex set. A simple example of
(\ref{eq:uncertainty_set}) for $n=m=d=n_q=1$ and no state dependency is
illustrated in Figure~\ref{fig:uncertainty_set_example}, where:
\begin{equation*}
  W = 1,\,\,
  R =
  \begin{bmatrix}
    1 \\
    -1
  \end{bmatrix},\,\,
  r =
  \begin{bmatrix}
    1 \\ 1
  \end{bmatrix},\,\,
  L_1 = 1,\,\,
  \phi_1(|u_k|) = |u_k|.
\end{equation*}

Summarizing, the control problem is to chose an input $u_k$ at each time
$k\in\integers_+$ such that given $x_0\in\mathcal X$ and the system dynamics
(\ref{eq:DLTV}), the constraints (\ref{eq:state_constraint}) and
(\ref{eq:input_constraint}) are satisfied subject to any uncertainty specified
by (\ref{eq:uncertainty_constraint}). We assume that the control objective to be
minimized can be expressed as a convex function of $x_k$ and $u_k$.


%% file: control_law.tex
\section{Control Law Description}
\label{sec:control_law}

In this section we develop an RMPC law with a planning horizon of length
$N\in\integers_{++}$ that solves the control problem presented in
Section~\ref{sec:problem_formulation}. We begin with a set theoretic motivation
for how the state constraint (\ref{eq:state_constraint}) may be robustly
satisfied given the input constraint (\ref{eq:input_constraint}) and the
uncertainty (\ref{eq:uncertainty_constraint}).

\begin{definition}
  A compact convex set $\mathcal I$ is feasible robust controlled
  invariant (fRCI) for system (\ref{eq:DLTV}) and constraints
  (\ref{eq:state_constraint}), (\ref{eq:input_constraint}) and
  (\ref{eq:uncertainty_constraint}) if it satisfies the following definition:
  \begin{equation*}
    \mathcal I = \{x\in\mathcal X : \exists u\in\mathcal U\text{ s.t. }
    Ax+Bu+Dp\in\mathcal I\,\,\forall p\in\mathcal P(x,u)\}.
  \end{equation*}
\end{definition}

By constraining $x_k\in\mathcal I\subseteq\mathcal X$, the state constraint can
thus be feasibly satisfied for all time as long as $x_0\in\mathcal I$. Numerous
set-based iterative methods are available for computing the maximal volume
$\mathcal I$ \cite{Blanchini2015,Kvasnica2015,Rakovic2006}. Approaches like
\cite{Rakovic2006} that attempt to deal with state and input dependent
uncertainty, however, scale poorly beyond $n=2$ because they require computing
set differences and projections of polytopes. Both operations have
$\mathcal O(2^n)$ complexity \cite{Baotic2009}. For this reason we use a more
computationally efficient algorithm developed for independent uncertainty and
which avoids these two operations \cite{Kvasnica2015}. This requires using the
convex hull of our uncertainty model:
\begin{equation}
  \label{eq:uncertainty_constraint_convex_hull}
  p_k\in\Co\{\mathcal P(x_k,u_k) : x_k\in\mathcal X,\, u_k\in\mathcal U\}.
\end{equation}

While this introduces conservatism into the computation of $\mathcal I$ by
considering a larger uncertainty set, the approach is more scalable to higher
dimensions given existing methods for computing fRCI sets. Furthermore,
Corollary~\ref{cor:rci_check} in Section~\ref{sec:main_results} presents a
computationally efficient method for checking if $\mathcal X$ itself is fRCI, at
which point the computation of $\mathcal I$ can be avoided altogether as
described in Algorithm~\ref{alg:mpc}. The proof of
Theorem~\ref{thm:recursive_feasibility} in Section~\ref{sec:main_results} shows
that $\mathcal I$ is a convex set and may therefore be inner-approximated by a
polytope. We henceforth assume that this polytopic description is available:
\begin{equation}
  \label{eq:rci_set}
  \mathcal I \triangleq \{x\in\reals^n : Gx\le g\},
\end{equation}
where $G\in\reals^{n_g\times n}$ and $g\in\reals^{n_g}$.

Since $\mathcal I$ is fRCI, for each $x_k\in\mathcal I$
$\exists u_k\in\mathcal U$ such that $x_{k+1}\in\mathcal I$
$\forall p_k\in\mathcal P(x_k,u_k)$. We compute such a $u_k$ via a tightened
state constraint in an optimization problem. Letting $x_k$ be the current state,
the system (\ref{eq:DLTV}) has the following impulse response over $t=1:N$
future time steps:
\begin{equation}
  \label{eq:impulse_response}
  x_{k+t} = A^tx_k+\sum_{i=0}^{t-1}B^A_iu_{k+i}+\sum_{i=0}^{t-1}D^A_ip_{k+i},
\end{equation}
where $B^A_i\triangleq A^{t-1-i}B$ and $D^A_i\triangleq A^{t-1-i}D$ (the index
$t$ is omitted for notational simplicity). An invariance-preserving input
sequence over the $N$-step planning horizon satisfies:
\begin{equation}
  \label{eq:rpi_v1}
  G\left(
    A^tx_k+\sum_{i=0}^{t-1}B^A_iu_{k+i}+\sum_{i=0}^{t-1}D^A_ip_{k+i}
  \right)\le g,
\end{equation}
for all sequences $p_{k+i}\in\mathcal P(x_{k+i},u_{k+i})$, $i=0:t-1$ and
$t=1:N$. Similarly to \cite{Blanchini1990b,Kolmanovsky1998} we reformulate this
requirement by maximizing the left hand side of (\ref{eq:rpi_v1}):
\begin{equation}
  \label{eq:reformulation_of_frci_condition}
  G_j^\transp\bar x_{k+t}+\sum_{i=0}^{t-1}\max_{p_{k+i}\in\mathcal P(x_{k+i},u_{k+i})}G_j^\transp D^A_ip_{k+i}\le g_j,
\end{equation}
for $j=1:n_g$ and $t=1:N$, where $\bar x_{k+t}$ denotes the nominal state after $t$
time steps:
\begin{equation}
  \label{eq:nominal_state}
  \bar x_{k+t} \triangleq A^tx_k+\sum_{i=0}^{t-1}B^A_iu_{k+i}.
\end{equation}
Note that the maximization term in (\ref{eq:reformulation_of_frci_condition}) is
the support function $\sigma_{D^A_i\mathcal P(x_{k+i},u_{k+i})}(G_j)\in\reals_{+}$ where
$G_j$ is $\mathcal I$'s $j$-th facet's normal. This support function induces
constraint tightening.

A complication arises in evaluating the support function due to an algebraic
loop in the state dependent uncertainty. The set $\mathcal P(x_{k+i},u_{k+i})$
depends on $x_{k+i}$ which itself depends on the uncertainty sequence over time
steps $0:i-1$ in the planning horizon. Thus, $p_{k+i}$ becomes dependent on
$p_k,\dots,p_{k+i-1}$ and the maximum value that it can take involves the
maximization of a convex function over a non-convex domain. Simulation
experience with the example in Section~\ref{sec:example} shows that a convex
upper bound to this maximization is too conservative, thus we prefer to simplify
by considering instead the uncertainty set $\mathcal P(\bar x_{k+i},u_{k+i})$
based on the nominal state. While this choice has no impact at $t=1$ (since
$x_k=\bar x_k$), the same cannot be said over the remaining $N-1$ steps because
perturbed states with a larger norm may occur, inducing a larger uncertainty
than the one predicted by $\mathcal P(\bar x_{k+i},u_{k+i})$. However, since
RMPC is implemented in receding horizon fashion and $u_k$ is anyway the only
input to be applied, the implementation remains robust. Furthermore, note that
the formulation remains exact over the entire planning horizon for input
dependent uncertainty.
With this in mind, the support function in
\eqref{eq:reformulation_of_frci_condition} is simplified using the separable
nature of (\ref{eq:uncertainty_set}) and H\"older's inequality:
\begin{equation}
  \label{eq:support_function_simplification}
  \begin{split}
    &\sigma_{D^A_i\mathcal P(\bar x_{k+i},u_{k+i})}(G_j)
    = \sigma_{D^A_i\mathcal W}(G_j)+\\
     &\sum_{l=1}^{n_q}\|G_j^\transp D^A_iL_l\|_{q_{q,l}}\phi_l(\|F_{x,l}\bar x_{k+i}\|_{p_{x,l}},\|F_{u,l}u_{k+i}\|_{p_{u,l}}),
   \end{split}
   \hspace{-1mm}
\end{equation}
where the $q_{q,l}$-norm is dual to the $p_{q,l}$-norm in
(\ref{eq:uncertainty_set}), i.e. $1/q_{q,l}+1/p_{q,l}=1$, and we denote by
$\mathcal W\triangleq\{Ww : Rw\le r\}$ the polytopic independent uncertainty
part of $\mathcal P(\bar x_{k+i},u_{k+i})$. Note that
(\ref{eq:support_function_simplification}) is not conservative since H\"older's
inequality is tight \cite{Boyd2004}. Using
(\ref{eq:support_function_simplification}), we can write
(\ref{eq:reformulation_of_frci_condition}) as a set of $Nn_g$ convex
constraints:
\begin{equation}
  \label{eq:rpi_v3}
  \begin{split}
    &G_j^\transp\bar x_{k+t}+
    \sum_{i=0}^{t-1}\sigma_{D^A_i\mathcal W}(G_j)+ \\
    &\sum_{i=0}^{t-1}\sum_{l=1}^{n_q}\|G_j^\transp D^A_iL_l\|_{q_{q,l}}\phi_l(\|F_{x,l}\bar x_{k+i}\|_{p_{x,l}}, \\
    &\hspace{37.5mm}\|F_{u,l}u_{k+i}\|_{q_{u_l}})
    \le g_j,
  \end{split}
  \hspace{-1mm}
\end{equation}
for $j=1:n_g$ and $t=1:N$. Importantly, $\sigma_{D^A_i\mathcal W}(G_j)$ can be
pre-computed offline leading to a more efficient implementation. The RMPC
on-line optimization problem to be implemented in receding horizon fashion can
thus be written as:
\begin{equation}
  \label{eq:mpc_problem}
  \minoptconstrained{
    u_{k},\dots,u_{k+N-1}
  }{
    J(\bar x_{k+1},\dots,\bar x_{k+N},u_k,\dots,u_{k+N-1})
  }{
    \text{(\ref{eq:input_constraint}) and (\ref{eq:rpi_v3}).}
  }
\end{equation}
where the cost function $J:(\reals^n)^N\times(\reals^m)^N\to\reals$ is a design
choice which should be convex. If $J$ is linear and the 1- or $\infty$-norms are
used for all $\phi_l$ in (\ref{eq:rpi_v3}), this is a linear program (LP). If
$J$ is quadratic, it is a QP. In all cases when the 2-norm is used in the
$\phi_l$ functions, (\ref{eq:mpc_problem}) is an SOCP. Note that all of these
problem types are convex and have efficient solvers available that guarantee
convergence to the global optimum when the feasible set is non-empty (which it
is certified to always be in Section~\ref{sec:main_results}). Furthermore, the
constraint count of (\ref{eq:mpc_problem}) is $\mathcal O(N)$ owing to the
welcome property that the effect of worst-case uncertainty on a discrete time
linear system is explicitly given by (\ref{eq:support_function_simplification}),
avoiding a combinatorial search. Algorithm~\ref{alg:mpc} summarizes the off-line
and on-line steps that make up the full RMPC controller.

\begin{algorithm}[t]
  \centering
  \begin{algorithmic}
    \State {\color{gray}\textbf{Off-line}:}
    \If{Corollary~\ref{cor:rci_check} is satisfied for $\mathcal I=\mathcal X$}
    \State $\mathcal I\gets\mathcal X$
    \Else
    \State Compute $\mathcal I$ using the set-based algorithm in \cite{Kvasnica2015}
    \If{Corollary~\ref{cor:rci_check} is not satisfied for $\mathcal I$}
    \State Report error ``$N$ is too large''
    \EndIf
    \EndIf
    \State Store $\sigma_{D^A_i\mathcal W}(G_j)$ for $j=1:n_g$,
    $t=1:N$, $i=0:t-1$
    \State {\color{gray}\textbf{On-line}:}
    \State Obtain the current state $x$
    \State Set $x_k\gets x$, solve (\ref{eq:mpc_problem}) and apply $u_k$
    \State Sleep $T_{\text{s}}$ seconds\Comment{Discretization time step}
    \caption{Off-line and on-line RMPC steps.}
    \label{alg:mpc}
  \end{algorithmic}
\end{algorithm}


%% file: main_results.tex
\section{Main Results}
\label{sec:main_results}

In this section, Theorem~\ref{thm:recursive_feasibility} proves that the control
law (\ref{eq:mpc_problem}) is recursively feasible and
Corollary~\ref{cor:rci_check} provides an alternative sufficient condition for
certifying $\mathcal X$ to be fRCI. This helps to avoid laborious set-based
approaches for computing $\mathcal I$.

\begin{theorem}
  \label{thm:recursive_feasibility}
  The optimization problem (\ref{eq:mpc_problem}) is recursively feasible if and
  only if it is feasible at the vertices of $\mathcal I$.
  \begin{proof}
    Let $\mathcal V\triangleq\{v_1,...,v_M\}$ be the set of vertices of
    $\mathcal I$. The forward implication is trivial. If (\ref{eq:mpc_problem})
    is recursively feasible then it is feasible in particular when
    $x_k\in\mathcal V$. For the reverse implication, suppose that we solve
    (\ref{eq:mpc_problem}) with $x_k$ set to each vertex $v_m$, $m=1:M$, and
    obtain associated optimal input sequences $u_{k+t}^m$, $t=0:N-1$, and
    nominal state sequences $\bar x_{k+t}^m$, $t=0:N$, as given by
    (\ref{eq:nominal_state}):
    \begin{equation*}
      \bar x_{k+t}^m = A^tv_m+\sum_{i=0}^{t-1}B^A_iu_{k+i}^m.
    \end{equation*}
    
    Consider now a state $x_k\in\mathcal I$. Since $\mathcal I$ is a convex
    polytope, we can express $x_k$ as a convex combination of $\mathcal I$'s
    vertices:
    \begin{equation*}
      x_k=\sum_{m=1}^M\theta_mv_m,\,\sum_{m=1}^M\theta_m=1,\,\theta_m\ge 0\,\,\forall m=1:M.
    \end{equation*}

    It is now possible to sum the constraint (\ref{eq:rpi_v3}) applied at each
    vertex, weighted by $\theta_m$, to obtain:
    \begin{equation*}
      \begin{split}
        &G_j^\transp\sum_{m=1}^M\theta_m\bar x_{k+t}^m+
        \sum_{i=0}^{t-1}\sigma_{D^A_i\mathcal W}(G_j)+ \\
        &\sum_{i=0}^{t-1}\sum_{l=1}^{n_q}\|G_j^\transp D^A_iL_l\|_{q_{q,l}}\hspace{-1mm}\sum_{m=1}^M\theta_m\phi_l(\|F_{x,l}\bar x_{k+i}^m\|_{p_{x,l}}, \\
        &\hspace{48mm}\|F_{u,l}u_{k+i}^m\|_{p_{u,l}}) \le g_j.
      \end{split}
    \end{equation*}

    Because each $\phi_l$ is convex, it follows from Jensen's inequality that:
    \begin{equation*}
      \begin{split}
        \phi_l(\|F_{x,l}&\sum_{m=1}^M\theta_m\bar x_{k+i}^m\|_{p_{x,l}},\|F_{u,l}\sum_{m=1}^M\theta_mu_{k+i}^m\|_{p_{u,l}})\le \\
        &
        \sum_{m=1}^M\theta_m\phi_l(\|F_{x,l}\bar x_{k+i}^m\|_{p_{x,l}},\|F_{u,l}u_{k+i}^m\|_{p_{u,l}}),
      \end{split}
    \end{equation*}
    and as a result we have:
    \begin{equation}
      \label{eq:rci_constraint_rec_feas_proof}
      \begin{split}
        &G_j^\transp\sum_{m=1}^M\theta_m\bar x_{k+t}^m+
        \sum_{i=0}^{t-1}\sigma_{D^A_i\mathcal W}(G_j)+ \\
        &\sum_{i=0}^{t-1}\sum_{l=1}^{n_q}\|G_j^\transp D^A_iL_l\|_{q_{q,l}}\phi_l(\|F_{x,l}\sum_{m=1}^M\theta_m\bar x_{k+i}^m\|_{p_{x,l}}, \\
        &\|F_{u,l}\sum_{m=1}^M\theta_mu_{k+i}^m\|_{p_{u,l}}) \le g_j.
      \end{split}
    \end{equation}

    Because $\mathcal U$ is convex and $u_{k+i}^m\in\mathcal U$ $\forall m=1:M$,
    $i=0:N-1$, we have $\sum_{m=1}^M\theta_mu_{k+i}^m\in\mathcal U$
    $\forall i=0:N-1$. Furthermore:
    \begin{equation*}
      \sum_{m=1}^M\theta_m\bar x_{k+t}^m = A^tx_k+\sum_{i=0}^{t-1}B^A_i\sum_{m=1}^M\theta_mu_{k+i}^m.
    \end{equation*}
    We thus recognize that (\ref{eq:rci_constraint_rec_feas_proof}) is nothing
    but (\ref{eq:rpi_v3}) for the initial state $x_k$ and the feasible input
    sequence $\sum_{m=1}^M\theta_mu_{k+i}^m$, $i=0:N-1$. The feasible set of
    (\ref{eq:mpc_problem}) is thus non-empty. Since this control input ensures
    that $x_{k+1}\in\mathcal I$ under the worst-case disturbance,
    (\ref{eq:mpc_problem}) continues to be feasible at the next time step which
    means that it is recursively feasible.
  \end{proof}
\end{theorem}

The following corollary provides a more computationally efficient method than
computational geometry approaches for verifying that a particular polytope is
fRCI. This is used in Algorithm~\ref{alg:mpc} for checking if $\mathcal X$ in
(\ref{eq:state_constraint}) is fRCI, avoiding an unnecessary call to the more
time consuming set-based algorithm in \cite{Kvasnica2015}, especially for high
dimensional systems.

\begin{corollary}
  \label{cor:rci_check}
  Let $\mathcal I=\mathcal X$. A sufficient condition for $\mathcal I$ to be
  fRCI is for (\ref{eq:mpc_problem}) to be feasible at its vertices. The
  condition becomes also necessary when $N=1$.
  \begin{proof}
    The sufficient condition follows from
    Theorem~\ref{thm:recursive_feasibility} and the fact that (\ref{eq:rpi_v3})
    yields the fRCI property by design whenever (\ref{eq:mpc_problem}) is
    recursively feasible. The necessity of this condition when $N=1$ is a
    consequence of that no conservatism in the planning problem is then present.
  \end{proof}
\end{corollary}


%% file: extensions.tex
\section{Extension to Semi-feedback RMPC}
\label{sec:extensions}

The open-loop RMPC law \eqref{eq:mpc_problem} can be overly conservative because
it does not model feedback action, meaning that the open-loop control policy
$\{u_k,\dots,u_{k+N-1}\}$ cannot mitigate the effect of disturbances during the
planning stage \cite{Bemporad2007}. Indeed, the disturbance action terms
$\sigma_{D_i^A\mathcal W}(G_j)$ and $\|G_j^\transp D_i^AL_l\|_{q_{q,l}}$ in
\eqref{eq:rpi_v3} are independent of the control policy.
This section extends our formulation to semi-feedback RMPC
\cite{Bemporad1998,Cannon2005} which introduces feedback action into the
planning stage without significantly increasing computational complexity. To
this end, a static feedback gain $K$ is designed and the control is re-defined
as
\begin{equation}
  \label{eq:semifeedback_control}
  u_k = v_k+Kx_k,
\end{equation}
such that $v_k\in\reals^m$ becomes the decision variable. The nominal state
\eqref{eq:nominal_state} becomes:
\begin{equation}
  \label{eq:sf_nominal_state}
  \bar x_{k+t} = \tilde A^tx_k+\sum_{i=0}^{t-1}B_i^{\tilde A}v_{k+i},
\end{equation}
where $\tilde A=A+BK$. The robust constraint \eqref{eq:rpi_v3} becomes:
\begin{equation}
  \label{eq:sf_rpi_v3}
  \begin{split}
    &G_j^\transp\bar x_{k+t}+
    \sum_{i=0}^{t-1}\sigma_{D^{\tilde A}_i\mathcal W}(G_j)+ \\
    &\sum_{i=0}^{t-1}\sum_{l=1}^{n_q}\|G_j^\transp D^{\tilde A}_iL_l\|_{q_{q,l}}\phi_l(\|F_{x,l}\bar x_{k+i}\|_{p_{x,l}}, \\
    &\hspace{27.5mm}\|F_{u,l}(v_{k+i}+K\bar x_{k+i})\|_{q_{u_l}})
    \le g_j,
  \end{split}
  \hspace{-1mm}
\end{equation}
where we have used the nominal state for the input-dependent input
error. Applying the discussion in Section~\ref{sec:control_law}, this choice has
no impact at $t=1$ so the RMPC law remains robust when implemented in receding
horizon fashion. With these changes, semi-feedback RMPC can be implemented in
the same way as \eqref{eq:mpc_problem}.


%% file: example.tex
\section{Illustrative Example}
\label{sec:example}

\begin{figure}
  \centering
  \includegraphics[width=0.4\columnwidth]{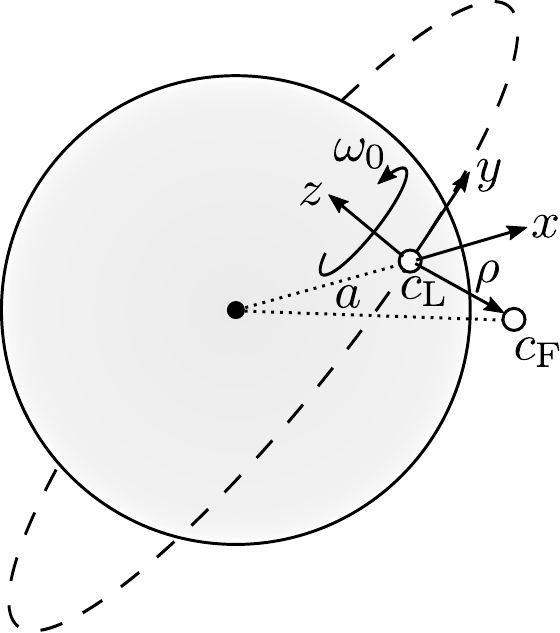}
  \vspace{\captionskip}
  \caption{LVLH frame showing the leader and follower satellites.}
  \label{fig:lvlh}
\end{figure}

In this section the control law presented in Sections~\ref{sec:control_law} and
\ref{sec:extensions} is applied to satellite robust position control in low
Earth orbit. Consider a two satellite formation as shown in
Figure~\ref{fig:lvlh} with a leader $c_{\text{L}}$ in a circular orbit with
radius $a$ and a follower $c_{\text{F}}$. The follower's translation in a local
vertical local horizontal (LVLH) frame is given by the Clohessy-Wiltshire
equations (time is omitted for notational simplicity):
\begin{subequations}
  \begin{align}
    \ddot x &= 3\omega_0^2x+2\omega_0\dot y+u_x+w_x, \label{eq:cw_x} \\
    \ddot y &= -2\omega_0\dot x+u_y+w_y, \label{eq:cw_y} \\
    \ddot z &= -\omega_0^2z+u_z+w_z, \label{eq:cw_z}
  \end{align}
\end{subequations}
where $\omega_0=\sqrt{\mu/a^3}$ and $\mu$ is the standard gravitational
parameter. Let $\rho=(x,y,z)\in\reals^3$ denote the follower's position relative
to the leader, where the coordinate $x$ is not to be confused with the state
vector. Using the state $x\triangleq (\rho,\dot\rho)\in\reals^6$, input
$u=(u_x,u_y,u_z)\in\reals^3$ and exogenous disturbance
$w=(w_x,w_y,w_z)\in\reals^3$, the dynamics (\ref{eq:cw_x})-(\ref{eq:cw_z}) take
on the familiar linear form:
\begin{equation}
  \label{eq:satellite_continuous_time}
  \dot x = A_{\text{c}}x+B_{\text{c}}(u+w).
\end{equation}

We use an impulsive control input in which the RCS system can induce an
instantaneous velocity increment $\Delta v(\tau)$ at time $\tau$ via the control
input $u(t)=\Delta v(\tau)\delta(t-\tau)$, where $\delta$ is the Dirac delta,
every $T_{\text{s}}$ seconds. Meanwhile, $w$ is assumed to be a constant
acceleration over the $T_{\text{s}}$ period. These assumptions allow
(\ref{eq:satellite_continuous_time}) to be discretized \cite{Antsaklis2007}:
\begin{equation}
  \label{eq:satellite_dlti}
  x_{k+1} = Ax_k+Bu_k+Ew_k,
\end{equation}
\begin{equation*}
  A = e^{A_{\text{c}}T_{\text{s}}},\quad
  B = AB_{\text{c}}, \quad
  E = \int_{0}^{T_{\text{s}}}e^{A(T_{\text{s}}-t)}B\mathrm{d}t.
\end{equation*}

We include three uncertainty sources:
\begin{enumerate}
\item Atmospheric drag, modeled as an independent component
  $w_k\in\{w : \|w\|_\infty\le w_{\max}\}$;
\item Additive input error using the Gates model, which captures the error in
  the RCS system's reproduction of a $\Delta v$ desired velocity increment
  \cite{Gates1963,Wagner2014,Goodson2013}. Because (\ref{eq:uncertainty_set})
  does not capture directional dependency, we confine ourselves to an isotropic
  description. In fact this is anyway the best modeling choice if the
  satellite's design is unknown \cite{Gates1963}. The error term is
  $v_k=v_k^{\text{fix}}+v_k^{\text{prop}}$ where
  $v_k^{\text{fix}}\in\{v : \|v\|_2\le\sigma_{\text{fix}}\}$ and
  $v_k^{\text{prop}}\in\{v : \|v\|_2\le\sigma_{\text{rcs}}\|u_k\|_2\}$;
\item Additive state estimation error
  $e_k=e_k^{\text{fix}}+(I,0)e_k^{\text{pos}}+(0,I)e_k^{\text{vel}}$ where:
\end{enumerate}
\begin{equation*}
  \begin{split}
    e_k^{\text{fix}} & \setlength\arraycolsep{3pt}
    \in\{e\in\reals^6 : \|\begin{bmatrix} I & 0
    \end{bmatrix}e\|_\infty\le p_{\max},
    \|\begin{bmatrix}
      0 & I
    \end{bmatrix}e\|_\infty\le v_{\max} \}, \\
    e_k^{\text{pos}} &\setlength\arraycolsep{3pt}
    \in\{e\in\reals^3 : \|e\|_\infty\le\sigma_{\text{pos}}\|
    \begin{bmatrix}
      I & 0
    \end{bmatrix}
    x_k\|_2\}, \\
    e_k^{\text{vel}} &\setlength\arraycolsep{3pt}
    \in\{e\in\reals^3 : \|e\|_\infty\le\sigma_{\text{vel}}\|
    \begin{bmatrix}
      0 & I
    \end{bmatrix}
    x_k\|_2\}.
  \end{split}
\end{equation*}

Altogether, these form a set of independent and dependent polytopic and
ellipsoidal uncertainties that is readily described by
(\ref{eq:uncertainty_set}). In particular, $d=21$, $n_q=4$,
$w=(w_k,e_k^{\text{fix}})$, $q_1=v_k^{\text{fix}}$, $q_2=v_k^{\text{prop}}$,
$q_3=e_k^{\text{pos}}$, $q_4=e_k^{\text{vel}}$, $p_{q,1}=p_{q,2}=2$,
$p_{q,3}=p_{q,4}=\infty$, $\phi_1=\sigma_{\text{fix}}$,
$\phi_2:u_k\mapsto\sigma_{\text{rcs}}\|u_k\|_2$,
$\setlength\arraycolsep{3pt}\phi_3:x_k\mapsto\sigma_{\text{pos}}\|\begin{bmatrix} I & 0
\end{bmatrix}x_k\|_2$, $\setlength\arraycolsep{3pt}\phi_4:x_k\mapsto\sigma_{\text{vel}}\|\begin{bmatrix}
  0 & I
\end{bmatrix}x_k\|_2$ and the following matrices:
\begin{gather*}
  D =
  \begin{bmatrix}
    E & -A & B & B & -A & -A
  \end{bmatrix} \\
  W = (I_{9},0_{18\times 9})\,\,\,\,
  R = (I_{9},-I_{9}) \\
  r = (w_{\max}\bm{1},p_{\max}\bm{1},v_{\max}\bm{1},w_{\max}\bm{1},p_{\max}\bm{1},v_{\max}\bm{1}) \\
  L_1 = (0_{9\times 3},I,0_{15\times 3})\,\,\,\,
  L_2 = (0_{12\times 3},I,0_{12\times 3}) \\
  L_3 = (0_{15\times 3},I,0_{9\times 3})\,\,\,\,
  L_4 = (0_{24\times 3},I),
\end{gather*}
where $\bm{1}\in\reals^3$ is a vector of ones. We use the following cost:
\begin{equation}
  \label{eq:satellite_cost}
  J\triangleq \sum_{t=0}^{N-1}\hat u_{k+t}^\transp\hat u_{k+t}+
  \lambda\hat{\bar x}_{k+t+1}^\transp\hat{\bar x}_{k+t+1},
\end{equation}
where $\hat u_{k+t}$ and $\hat{\bar x}_{k+t+1}$ are the scaled input and nominal
state such that they attain plus or minus unity at the boundary of their
respective constraint polytopes $\mathcal U$ and $\mathcal I$ while
$\lambda=0.003$ is a manually chosen trade-off weight. The input penalty
reflects a minimum-fuel type problem and the state penalty endows the finite
horizon control law with an otherwise lacking long-term knowledge that the
origin corresponds to minimal fuel usage, since nominally it takes zero control
action to remain there. In this problem, $\mathcal X$ in
(\ref{eq:state_constraint}) takes on the direct interpretation of a maximum
control error specification. We use the following numerical values:
\begin{gather*}
  \setlength\arraycolsep{2pt}
  \mathcal X=\{x\in\reals^6 : \|\begin{bmatrix}
    I & 0
  \end{bmatrix}x\|_\infty\le 10\text{ cm},
  \|\begin{bmatrix}
    0 & I
  \end{bmatrix}x\|_\infty\le 1\text{ mm/s}\} \\
  \mathcal U = \{u\in\reals^3 : \|u\|_\infty\le 2\text{ mm/s}\} \\
  \mu = 3.986\cdot 10^{14}\text{ m}^3\text{/s}^2\,\,\,
  a = 6793.137\text{ km}\,\,\,
  T_{\text{s}} = 100\text{ s} \\
  w_{\max} = 50\text{ nm/s}^2\,\,\,
  \sigma_{\text{fix}} = 1\text{ }\mu\text{m/s}\,\,\,
  p_{\max} = 0.4\text{ cm}\,\,\,N=4 \\
  v_{\max} = 4\text{ }\mu\text{m/s}\,\,\,
  \sigma_{\text{rcs}} = \tan\frac{\pi}{180}\,\,\,
  \sigma_{\text{pos}} = 0.02\,\,\,
  \sigma_{\text{vel}} = 0.001.
\end{gather*}

We compare the following four controllers:
\begin{enumerate}
\item Nominal MPC: ignores the uncertainty, effectively removing the summations
  in (\ref{eq:rpi_v3});
\item Conservative RMPC: replaces $\phi_i$ in (\ref{eq:uncertainty_set}) by its
  maximum over $x_k\in\mathcal X$ and $u_k\in\mathcal U$, which results in a
  conservative independent uncertainty model;
\item Our open-loop RMPC law (\ref{eq:mpc_problem});
\item Our semi-feedback RMPC law using the modifications described in
  Section~\ref{sec:extensions}. We design $K$ via LQR with weight matrices
  $Q=I_{6}$ and $R=10^5I$, scaled to $\mathcal I$ and $\mathcal U$ respectively
  in the same way as for \eqref{eq:satellite_cost}.
\end{enumerate}

\begin{figure}
  \centering
  \includegraphics[width=1\columnwidth]{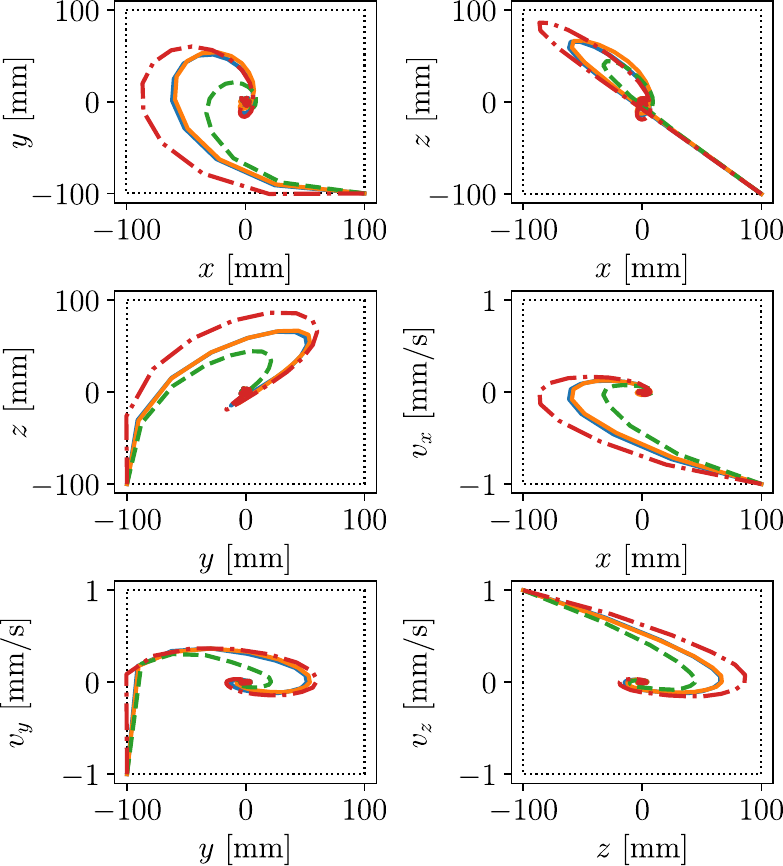}
  \vspace{\captionskip}
  \caption{Projections of the transient response shown in dash dotted red for
    nominal MPC, dashed green for conservative RMPC, solid blue for our
    open-loop RMPC and solid orange for our semi-feedback RMPC. The dotted black
    rectangle shows the boundary of $\mathcal I$.}
  \label{fig:transient_response}
\end{figure}

In each case the satellite is acted upon by all three of the uncertainty sources
described above. It turns out that for this problem,
Corollary~\ref{cor:rci_check} succeeds and so $\mathcal I=\mathcal X$. To
demonstrate the conservatism exhibited by the conservative RMPC law, consider
Figure~\ref{fig:transient_response} which shows the four controllers' transient
response when starting from a vertex of $\mathcal I$. While nominal MPC briefly
exits $\mathcal I$, conservative RMPC tends to drive the satellite more quickly
into the interior of $\mathcal I$ because it assumes more uncertainty than
necessary. Open-loop and semi-feedback RMPC yield similar responses that are
somewhere in between nominal and conservative MPC, staying within $\mathcal I$
and not avoiding the boundaries of $\mathcal I$ unnecessarily.

The RMPC law (\ref{eq:mpc_problem}) enables the control engineer to work with a
richer set of feasible parameters. For example, the required position accuracy
can be increased to 5~cm without changing any other parameter, while doing so
with the conservative RMPC law requires using $N\le 2$. 
The open-loop RMPC law works for $N\le 4$ before the passively propagated
uncertainty ``outgrows'' $\mathcal X$. The less conservative semi-feedback RMPC
law works for $N\le 6$. 

\begin{figure}
  \centering
  \includegraphics[width=1\columnwidth]{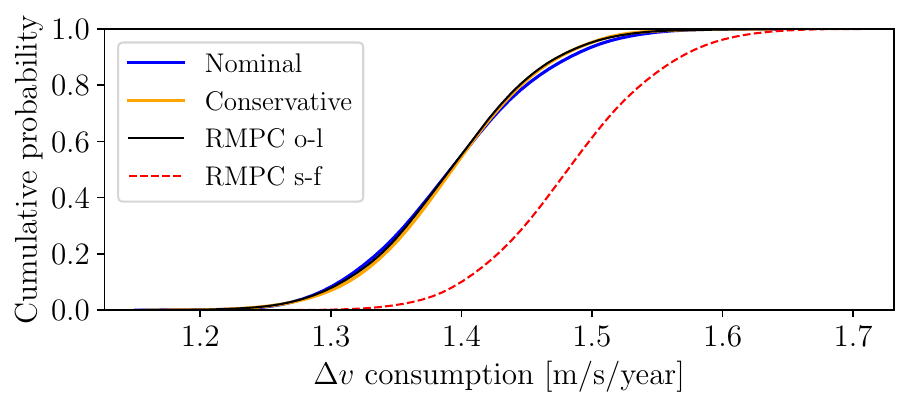}
  \vspace{\captionskip}
  \caption{Fuel usage distribution for each controller from 2000 Monte Carlo
    simulations. \textit{RMPC o-l} and \textit{s-f} stand for our open-loop and
    semi-feedback RMPC laws respectively.}
  \label{fig:deltavslope_distribution}
\end{figure}

\begin{figure}
  \centering
  \includegraphics[width=1\columnwidth]{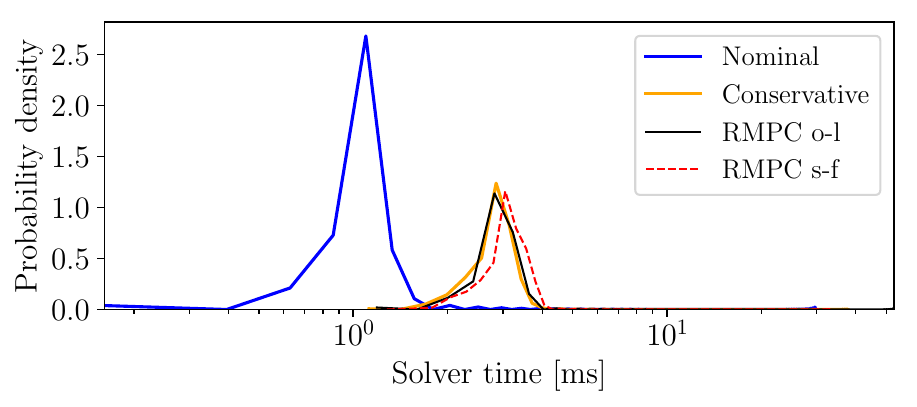}
  \vspace{\captionskip}
  \caption{Optimization solver time distribution for each controller from 2000
    Monte Carlo simulations. \textit{RMPC o-l} and \textit{s-f} stand for our
    open-loop and semi-feedback RMPC laws respectively.}
  \label{fig:solver_time}
\end{figure}

Next, we compare the fuel consumption and computational efficiency of the three
controllers. Fuel consumption is quantified by the sum of velocity increment
magnitudes commanded per year as obtained from a linear regression, yielding
units of m/s/year. This is common in space applications, where a mission is
characterized by a ``$\Delta v$ budget''. Computational efficiency is measured
by the time taken to solve (\ref{eq:mpc_problem}). Since both quantities are
affected by the uncertainty realization, we run 2000 Monte Carlo simulations for
each controller where the satellite is initialized at the origin and is
controlled for a duration of four
orbits.

Figure~\ref{fig:deltavslope_distribution} shows a fuel consumption cumulative
distribution plot using a Gaussian smoothing kernel. There is no statistically
significant difference between the nominal MPC, conservative RMPC and open-loop
RMPC laws. As expected, the embedded LQR controller increases the fuel
consumption of semi-feedback RMPC by an average amount of $0.09$~m/s/year.

Figure~\ref{fig:solver_time} shows the solver time distribution for Python
2.7.15 with ECOS \texttt{2.0.7.post1} \cite{Domahidi2013} in Ubuntu 18.04.1 with
a 3.60~GHz Intel i7-6850K CPU and 64~GB of RAM. All four controllers have
comparable solver times. The distributions reflect the problem difficulty
hierarchy, wherein our RMPC laws are the most difficult due to second order cone
constraints. On average, nominal MPC takes 1.1~ms, conservative and open-loop
RMPC take 2.9 ms, and semi-feedback RMPC takes 3~ms. These differences are all
statistically significant at the 99.9\% confidence level. Generally speaking,
SOCP problems are amenable to real-time implementation and we expect this to be
the case here \cite{Dueri2014,Dueri2017,Zeilinger2014,Domahidi2013}.



%% file: conclusion.tex
\section{Conclusion}
\label{sec:conclusion}

This paper introduced a state and input dependent uncertainty model and a
corresponding computationally efficient robust receding horizon control law
based on second order cone programming. We have shown that the control law is
recursively feasible and that both open-loop and semi-feedback implementations
are possible. Simulations of a satellite system demonstrate that the approach is
more versatile than RPMC based on an independent uncertainty model. The
approach, however, is somewhat hampered by the ability to compute a robust
controlled invariant set for the system. Recent work on controllable set
computation for convex optimal control problems has the potential to remove this
shortcoming \cite{Dueri2016}.


%% file: root.bbl
\begin{thebibliography}{10}

\bibitem{Gates1963}
C.~R. Gates, ``A simplified model of midcourse maneuver execution errors,''
  Tech. Rep. 32-504, JPL, Pasadena, CA, oct 1963.

\bibitem{Wagner2014}
S.~V. Wagner, ``Maneuver performance assessment of the {Cassini} spacecraft
  through execution-error modeling and analysis,'' in {\em AAS/AIAA Space
  Flight Mechanics Meeting}, jan 2014.

\bibitem{Navvabi2018}
H.~Navvabi and A.~H.~D. Markazi, ``New {AFSMC} method for nonlinear system with
  state-dependent uncertainty: Application to hexapod robot position control,''
  {\em Journal of Intelligent {\&} Robotic Systems}, may 2018.

\bibitem{Bemporad2007}
A.~Bemporad and M.~Morari, ``Robust model predictive control: A survey,'' in
  {\em Robustness in identification and control}, pp.~207--226, Springer
  London, 2007.

\bibitem{Mayne2014}
D.~Q. Mayne, ``Model predictive control: Recent developments and future
  promise,'' {\em Automatica}, vol.~50, pp.~2967--2986, dec 2014.

\bibitem{Lee1997}
J.~Lee and Z.~Yu, ``Worst-case formulations of model predictive control for
  systems with bounded parameters,'' {\em Automatica}, vol.~33, pp.~763--781,
  may 1997.

\bibitem{Bemporad1998}
A.~Bemporad, ``Reducing conservativeness in predictive control of constrained
  systems with disturbances,'' in {\em Proceedings of the 37th {IEEE}
  Conference on Decision and Control (Cat. No.98CH36171)}, {IEEE}, 1998.

\bibitem{Scokaert1998}
P.~Scokaert and D.~Q. Mayne, ``Min-max feedback model predictive control for
  constrained linear systems,'' {\em {IEEE} Transactions on Automatic Control},
  vol.~43, no.~8, pp.~1136--1142, 1998.

\bibitem{Kerrigan2004}
E.~C. Kerrigan and J.~M. Maciejowski, ``Feedback min-max model predictive
  control using a single linear program: robust stability and the explicit
  solution,'' {\em International Journal of Robust and Nonlinear Control},
  vol.~14, no.~4, pp.~395--413, 2004.

\bibitem{Lofberg2003}
J.~L\"ofberg, {\em {Minimax approaches to robust model predictive control}}.
\newblock {Dissertation (Ph.D.)}, Link\"oping University, 2003.

\bibitem{Blanchini1990b}
F.~Blanchini, ``Feedback control for linear time-invariant systems with state
  and control bounds in the presence of disturbances,'' {\em {IEEE}
  Transactions on Automatic Control}, vol.~35, no.~11, pp.~1231--1234, 1990.

\bibitem{Blanchini1990a}
F.~Blanchini, ``Control synthesis for discrete time systems with control and
  state bounds in the presence of disturbances,'' {\em Journal of Optimization
  Theory and Applications}, vol.~65, pp.~29--40, apr 1990.

\bibitem{Blanchini1990c}
F.~Blanchini, ``Constrained control for perturbed linear systems,'' in {\em
  29th {IEEE} Conference on Decision and Control}, {IEEE}, 1990.

\bibitem{Blanchini1994}
F.~Blanchini, ``Ultimate boundedness control for uncertain discrete-time
  systems via set-induced {Lyapunov} functions,'' {\em {IEEE} Transactions on
  Automatic Control}, vol.~39, no.~2, pp.~428--433, 1994.

\bibitem{Marruedo2002}
D.~Marruedo, T.~Alamo, and E.~Camacho, ``Input-to-state stable {MPC} for
  constrained discrete-time nonlinear systems with bounded additive
  uncertainties,'' in {\em Proceedings of the 41st {IEEE} Conference on
  Decision and Control, 2002.}, {IEEE}, dec 2002.

\bibitem{Pin2009}
G.~Pin, D.~Raimondo, L.~Magni, and T.~Parisini, ``Robust model predictive
  control of nonlinear systems with bounded and state-dependent
  uncertainties,'' {\em {IEEE} Transactions on Automatic Control}, vol.~54,
  pp.~1681--1687, jul 2009.

\bibitem{Lee2002}
Y.~Lee, B.~Kouvaritakis, and M.~Cannon, ``Constrained receding horizon
  predictive control for nonlinear systems,'' {\em Automatica}, vol.~38,
  pp.~2093--2102, dec 2002.

\bibitem{Limon2004}
D.~Limon, T.~Alamo, and E.~Camacho, ``Robust stability of min-max {MPC}
  controllers for nonlinear systems with bounded uncertainties,'' in {\em
  Proceedings of the 16th Mathematical Theory of Networks and Systems
  Conference.}, 2004.

\bibitem{Acikmese2010}
B.~\Acikmese, J.~M. {Carson III}, and D.~S. Bayard, ``A robust model predictive
  control algorithm for incrementally conic uncertain/nonlinear systems,'' {\em
  International Journal of Robust and Nonlinear Control}, vol.~21,
  pp.~563--590, jul 2010.

\bibitem{Rakovic2009}
S.~V. Rakovi{\'{c}}, ``Set theoretic methods in model predictive control,'' in
  {\em Nonlinear Model Predictive Control}, pp.~41--54, Springer Berlin
  Heidelberg, 2009.

\bibitem{Cannon2005}
M.~Cannon and B.~Kouvaritakis, ``Optimizing prediction dynamics for robust
  {MPC},'' {\em {IEEE} Transactions on Automatic Control}, vol.~50,
  pp.~1892--1897, nov 2005.

\bibitem{Goodson2013}
T.~D. Goodson, ``Execution-error modeling and analysis of the {GRAIL}
  spacecraft pair,'' in {\em AAS/AIAA Space Flight Mechanics Meeting}, aug
  2013.

\bibitem{Dueri2014}
D.~Dueri, J.~Zhang, and B.~\Acikmese, ``Automated custom code generation for
  embedded, real-time second order cone programming,'' {\em {IFAC} Proceedings
  Volumes}, vol.~47, no.~3, pp.~1605--1612, 2014.

\bibitem{Dueri2017}
D.~Dueri, B.~\Acikmese, D.~P. Scharf, and M.~W. Harris, ``Customized real-time
  interior-point methods for onboard powered-descent guidance,'' {\em Journal
  of Guidance, Control, and Dynamics}, vol.~40, pp.~197--212, feb 2017.

\bibitem{Zeilinger2014}
M.~N. Zeilinger, D.~M. Raimondo, A.~Domahidi, M.~Morari, and C.~N. Jones, ``On
  real-time robust model predictive control,'' {\em Automatica}, vol.~50,
  pp.~683--694, mar 2014.

\bibitem{Domahidi2013}
A.~Domahidi, E.~Chu, and S.~Boyd, ``{ECOS}: An {SOCP} solver for embedded
  systems,'' in {\em 2013 European Control Conference ({ECC})}, {IEEE}, jul
  2013.

\bibitem{Boyd2004}
S.~Boyd and L.~Vandenberghe, {\em Convex Optimization}.
\newblock Cambridge University Press, 2004.

\bibitem{Blanchini2015}
F.~Blanchini and S.~Miani, {\em Set-Theoretic Methods in Control}.
\newblock Springer International Publishing, 2015.

\bibitem{Kvasnica2015}
M.~Kvasnica, B.~Tak{\'{a}}cs, J.~Holaza, and D.~Ingole, ``Reachability analysis
  and control synthesis for uncertain linear systems in {MPT},'' {\em
  {IFAC}-{PapersOnLine}}, vol.~48, no.~14, pp.~302--307, 2015.

\bibitem{Rakovic2006}
S.~V. Rakovi{\'{c}}, E.~C. Kerrigan, D.~Q. Mayne, and J.~Lygeros,
  ``Reachability analysis of discrete-time systems with disturbances,'' {\em
  {IEEE} Transactions on Automatic Control}, vol.~51, pp.~546--561, apr 2006.

\bibitem{Baotic2009}
M.~Baoti\'c, ``Polytopic computations in constrained optimal control,'' {\em
  Automatika, Journal for Control, Measurement, Electronics, Computing and
  Communications}, vol.~50, pp.~119--134, 2009.

\bibitem{Kolmanovsky1998}
I.~Kolmanovsky and E.~G. Gilbert, ``Theory and computation of disturbance
  invariant sets for discrete-time linear systems,'' {\em Mathematical Problems
  in Engineering}, vol.~4, no.~4, pp.~317--367, 1998.

\bibitem{Antsaklis2007}
P.~J. Antsaklis and A.~N. Michel, {\em A Linear Systems Primer}.
\newblock Birkh\"{a}user Boston, 2007.

\bibitem{Acikmese2007}
B.~\Acikmese and S.~R. Ploen, ``Convex programming approach to powered descent
  guidance for mars landing,'' {\em Journal of Guidance, Control, and
  Dynamics}, vol.~30, pp.~1353--1366, sep 2007.

\bibitem{Acikmese2011}
B.~\Acikmese and L.~Blackmore, ``Lossless convexification of a class of optimal
  control problems with non-convex control constraints,'' {\em Automatica},
  vol.~47, pp.~341--347, feb 2011.

\bibitem{Blackmore2012}
L.~Blackmore, B.~\Acikmese, and J.~M. {Carson III}, ``Lossless convexification
  of control constraints for a class of nonlinear optimal control problems,''
  {\em Systems {\&} Control Letters}, vol.~61, pp.~863--870, aug 2012.

\bibitem{Schouwenaars2006}
T.~Schouwenaars, {\em {Safe trajectory planning of autonomous vehicles}}.
\newblock {Dissertation (Ph.D.)}, Massachusetts Institute of Technology, 2006.

\bibitem{Dueri2016}
D.~Dueri, S.~V. Rakovi{\'{c}}, and B.~\Acikmese, ``Consistently improving
  approximations for constrained controllability and reachability,'' in {\em
  2016 European Control Conference ({ECC})}, {IEEE}, jun 2016.

\end{thebibliography}
